\documentclass[12pt,twoside,leqno]{article}
\usepackage{amssymb}
\newtheorem{Def}{Definition}[section]
\newtheorem{Pro}{Proposition}[section]
\newtheorem{Teo}{Theorem}[section]
\newtheorem{Lem}{Lemma}[section]
\newtheorem{Cor}{Corollary}[section]
\newtheorem{Rem}{Remark}[section]
\newenvironment{proof}[1][Proof]{\textbf{#1.} }{\ \hfill{\fbox{}}}
\newcommand{\A}{{\mathbb A}}
\newcommand{\C}{{\mathbb C}}

\newcommand{\I}{{\mathbb I}}
\newcommand{\bL}{{\mathbb L}}
\newcommand{\N}{{\mathbb N}}

\newcommand{\Z}{{\mathbb Z}}

\newcommand{\Sn}{{\mathcal S}_n}
\newcommand{\Snp}{{\mathcal S}_{n+1}}
\newcommand{\ds}{\displaystyle}
\begin{document}
\pagestyle{myheadings}
\markboth{\underline{\centerline{\textit{\small{O. EL
FOURCHI}}}}}{\underline{\centerline{\textit{\small{Spherical
analysis on permutations group and applications }}}}}
\title{\hspace{5cm} \underline{Spherical analysis on permutations group and applications
}}
\author{By\\ Omar El FOURCHI and Adil ECHCHELH\\Laboratoire de recherche en Management
 et m\'ethode quantitatives\\Facult\'e des Sciences Juridiques \'economiques et sociales\\K\'enitra-Morocco\\E.mail:
elfourchi$\_$omar$@$yahoo.com and echeladil@yahoo.fr}
\date{}
\maketitle

{\bf Summary}. The principal aim of the present paper is to
develop the theory of Gelfand pairs on the symmetric group in
order to define and study the horocyclic Radon transform on this
group. We also find a simple inversion formula for the Radon
transform of the solution to the heat equation associated to this
group.
\section{Introduction and preliminaries}
Radon transform play a critical role in subjects as diverse as
application to partial differential equations, X-ray technology
and radioastronomy. Like much of mathematics in the field of
harmonic analysis and integral geometry on homogeneous space has
some of its application in the work of Helgason cf.[12]. Although
the permutation group form one of the oldest parts of group theory
and  the harmonic analysis in this work may be regarded as a
trivial of compact case. Through the ubiquity of group actions and
representations theory, permutations group continue to be lively
topic of research in their own right (see [14],[17,18]). Working
on a single spherical analysis and developing the theory of
Gelfand pairs on permutations group we essentially want to define
the horocyclic Radon and construct its inverse. As an example of
the use Radon transforms, we give a simple solution to the Radon
transform of the heat equation on permutations group, all done
very explicitly.\\ A bijective function from
$Z_n=\{1,2,3,....,n\}$ onto itself is called a permutation of $n$
numbers; the set of all permutations of $n$ numbers, together with
the usual composition of functions,
 is called the symmetric group of degree $n$.
This group will be denoted by $\Sn$. Note that $\Sn$ is defined
for $n\geq 0$, and $\Sn$ has $n!$ elements (where $0! = 1$). If
$Y$ is a subset of $Z_n$, we shall write ${\mathcal{S}}_{Y}$ for
the subgroup of $\Sn$ which fixes every number outside $Y$.\\A
permutation $\sigma \in \Sn$, which interchanges two
 distinct numbers $i$ and $j$ and leaves all other numbers fixed,
 is called a transposition and is written as $\sigma = \tau_{i,j}$.
The function $\epsilon : \Sn \to \{\pm 1\}$,
 such that $\epsilon(\sigma)=(-1)^{N}$ if $\sigma$ is a product of $N$
 transpositions, is well-defined, called the signature of $\sigma$.
The number $\epsilon(\sigma)$ depends only on the parity of $N$
 and we have $\epsilon(\sigma\cdot\zeta)=\epsilon(\sigma).\epsilon(\zeta)$.
 \\The normalized Haar measure in $\Sn$ is given by
 $\nu=\frac{1}{n!}\sum_{\sigma \in \Sn}\delta_{\sigma}$,
 where $\delta_{\sigma}$ is Dirac measure.\\
The complex group algebra of the group $\Sn$ is $
\C(\Sn)=\left\{\sum_{\sigma \in \Sn}
 \lambda_{\sigma}
 \delta_{\sigma} :\lambda_{\sigma}\in \C\right\}$.\\
This is a vector space over $\C$, for which the set
$\{\delta_{\sigma}:\sigma \in \Sn\}$
 is a basis.
We note that the algebras $\C(\Sn)$,
 $\bL^{1}(\Sn)$ and $\bL^{2}(\Sn)$ (the space of integrable functions resp the square integrable functions )
 are all equal (see [4] and [5]).\\
The structure of this paper is the following. In Section 2, We
recall briefly the main definitions and results of the
representation theory on the symmetric group.\\ In Section 3, we
give a characterization of the set right cosets $\Snp \slash \Sn$
and those of double cosets $\Sn\backslash\Snp \slash \Sn $. The main
goal of this section will be devoted to spherical function of the
Gelfand pair $(\Snp, \Sn)$. In Section $4$ and $5$ we introduce and
invert the spherical Fourier transform in permuatation group. In
Section $6$ we introduce and investigate the Radon transform on
permutations group. In Section $7$ we establish the connection
between this transform and the solution to the heat equations
associated to this group which is the technical heart of the paper.

\section{Irreducible representations of permutations group}

We begin this section by recalling a few facts from the
representation theory of permutation groups cf. [14],[5]and
[17].\\ A subset $S$ of $Z_{n}$ is invariant under a permutation
$\sigma\in\Sn$
 if $\sigma(S)=S$.
A permutation is said to be circular when it admits
 $\O$ and $Z_{n}$ as its only invariant subsets.\\
An invariant subset $S$ of $Z_{n}$ is called a cycle of $\sigma$
 if $\sigma|_S$ is circular,
 where we are writing $\sigma \mid_{S}$ for the
 permutation which coincides with $\sigma$ on $S$,
 and which is the identity outside of $S$.

\begin{Lem}
Let  $\sigma \in \Sn$, the cycles of $\sigma$, denoted $S_1 , S_2
,..., S_k$,
 form a partition of $Z_{n}$.
Furthermore for all $i,j$,
 $\sigma\mid_{S_i}$ and $\sigma \mid_{S_j}$ commute and we have
$$
 \sigma
 =
 \prod_{i=1}^{k}
 \sigma \mid_{S_i}.
$$
\end{Lem}

Since the conjugacy class of an element $\sigma\in\Sn$ is
characterized by
 the lengths of the cycles of $\sigma$ (with repetitions),
 the number of conjugacy classes in $\Sn$ is equal to the number of
 partitions of $n$.
As (see [4]) the number of inequivalent
 complex irreducible representations of $\Sn$
 is equal to the number of conjugacy classes of $\Sn$.
Therefore the number of inequivalent complex irreducible
representations
 of $\Sn$ is equal to the number of partitions of $n$.

We should therefore aim to construct a regular representation of
 $\Sn$ for each partition of $n$.
This is made easier by knowing the primitive idempotents (see
[4]).

A non-zero element $e$ of $\C(\Sn)$ is said idempotent if $e*e=e$.
More generally $e$ is said essentially idempotent
 if $e*e=\lambda e$ for some $\lambda \neq 0$.
An idempotent $e$ is said to be primitive if $e$ decomposes
 as the sum of two idempotents: $e= e' + e''$ with
 $e'* e'' = e'' * e' = 0$.
An idempotent which does not decompose in this way is called a
 primitive idempotent.

In order to describe the primitive idempotents of $\C(\Sn)$,
 we will need to recall some definitions.\\
If $\lambda =(n_1,n_2,n_3,\ldots,n_k)$ is a partition of $n$
 (i.e., $n_1 \geq n_2 \geq n_3 \geq ....\geq n_k \geq 1$
 and $n = n_1 + n_2 +n_3 +...+n_k$),
 we associates to this partition a tableau
 $[\lambda ]=\{(i,j): i,j \in \Z; 1\leq i; 1\leq j \leq n_i \}$
 (here $\Z$ denotes the set of integers).\\
If $(i,j)\in [\lambda]$, then $(i,j)$ is called a node of
$[\lambda]$. The $k^{th}$ row (respectively column) of tableau
consists of those
 nodes whose first (respectively second) coordinate is $k$.\\
A Young-diagram is the one of the $n!$ arrays of integers obtained
by
 replacing each node in $[\lambda]$ by one of the integers
 $1,2,3,4,....,n$ allowing no repeats.\\
To the Young-diagram $t$, we associate its row-stabilizer,
 ${\mathcal{P}}_t$, is the subgroup of $\Sn$ keeping the rows of $t$
 fixed setwise.
i.e., $$
 {\mathcal{P}}_t
 =
 \{\sigma \in \Sn :
 \hbox{for all $i\in Z_n$, $i$ and $\sigma(i)$
 belong to the same row of $t$}\}.
$$ The column-stabilizer ${\mathcal{Q}}_t$, of $t$ is defined
similarly.
 \begin{Pro}(see [5])
Let $t$ a $\lambda$- tableau, if $\pi \in \Sn$ then
${\mathcal{P}}_{\pi t}=\pi {\mathcal{P}}_t \pi^{-1}$
${\mathcal{Q}}_{\pi t}=\pi {\mathcal{Q}}_t \pi^{-1}$
 \end{Pro}
 We define a relation of equivalence in the set of
 $\lambda$-tableau by\\
 $t_1 \sim t_2$ if and only if there exists $\pi \in
 {\mathcal{P}}_{t_1}$ such that $\pi t_1 =t_2$.\\
 The conjugacy class of tableau modulo this relation of
 equivalence are called tabloids and the conjugacy class of
 tableau $t$ is the tabloid noted by $\{t\}$.
\begin{Pro}(see [5])
The permutations group acts on the set of $\lambda$-tableau in the
following way: If $\pi \in \Sn$ and $t$ a $\lambda$-tableau  then
$\pi \{t\}= \{\pi t\}$
 \end{Pro}
 To each partition $\mu$ of $n$, $\mu = (p_1, p_2, ....,p_k)$, we
 associate the Young sub-group $S_\mu$ of $\Sn$ defined as product
 of the following sub-group
 $$S_{p_1+...p_{i-1}+1,p_1+...p_{i-1}+2,...,p_1+p_2+...+p_{i-1}+p_i},\,\,\,\,\,i=1,2,3,...,k.$$
we have then $$S_{1,...p_{1}}\times S_{p_1+1,...p_{1}+p_2,}\times
S _{p_1+p_2 +1,...,p_{1}+p_2+p_3}\times ....$$ If $\mu$ is a
partition of $n$, notice $M^{\mu}$ the ${\C}$vector space whose
basis are the $\mu$ distinct tabloids, so $M^{\mu}$ is $\C
[\Sn]$-module\\ For any Young-diagram $t$, we associate the
element of
 $\C(\Sn)$ defined as follows
$$
 e_t
 =
 \sum_{q\in{\mathcal{Q}}_t}
 \sum_{p\in {\mathcal{P}}_t}
 \epsilon(q)\delta_{p}*\delta_{q}.
$$ We note (see [4]), that $e_t$ is essentially idempotent and
 $\frac{1}{\lambda_t}e_t$ is a primitive idempotent ($\lambda_t \neq
 0$).
 In some notation $e_t$ is called polytabloid
 \begin{Lem}
If $\pi \in \Sn$ and if $e_t$ a polytabloid then $\pi e_t =e_{\pi
t }$,
\begin{Def}
The specht module $S^\mu$ is $\C[\Sn]$-module monogene generated
by any $\mu$-tabloids.\\ {\bf Remark.}Every result interpreted via
the Specht module is the same via the left ideal ${\mathcal{O}}_t$
of $\C(\Sn)$ generated by $e_t$ cf [14,p 17]
\end{Def}
\end{Lem}
\begin{Teo}(see [4] p :67) Let $t$ be a
Young-diagram, ${\mathcal{P}}_t$ its row-stabilizer,
 ${\mathcal{Q}}_t$ its column-stabilizer, and let $e_t$ be the element of
 $\C(\Sn)$ defined by
 $$
  e_t
  =
  \sum_{q\in {\mathcal{Q}}_t}
  \sum_{p\in{\mathcal{P}}_t}
  \epsilon(q)\delta_{p}*\delta_{q}.
 $$
We shall write ${\mathcal{O}}_t$ for the left ideal of $\C(\Sn)$
 generated by $e_t$ and ${\mathcal{R}}_t$ for the
 associated representation of $\Sn$.
Then
\begin{itemize}
\item
 ${\mathcal{R}}_t$ is irreducible;
\item
 two such representations ${\mathcal{R}}_s$ and ${\mathcal{R}}_t$
 are equivalent if and only if $s$ and $t$
 are Young diagrams for the same partition $\lambda$.
\end{itemize}
As the number of partitions is equal to the number of irreducible
 complex representations, we may obtain a representative ${\mathcal{R}}_t$
 for each equivalence class of irreducible representations by
 choosing for each partition $\lambda$ a Young diagram $t$.
\end{Teo}
\begin{Rem}
F. Scarabotti cf.[17] has given a short proof of a
characterization of James [14]of the irreducible modules as the
intersection of kernels of certain invariant operators using the
class sum of transpositions and a collection of related transform
for the complex representation of the permutations group.\end{Rem}
\section{Harmonic analysis of the pair $(\Snp ,\Sn)$}
We may regard $\Sn={\mathcal{S}}(Z_{n})$ as a subgroup of
 $\Snp ={\mathcal{S}}(Z_{n+1})$.
More precisely for $\sigma \in \Sn$, the map
 $\overline{\sigma}:Z_{n+1}\to Z_{n+1}$, defined by $x\mapsto \sigma(x)$ for $x\in Z_{n}$
 and $n+1 \mapsto n+1$, is an element of $\Snp$. We note that $\Sn$ acts transitively on $Z_{n}$ via the map
 $\Sn\times Z_{n} \to Z_{n}$: $(\sigma ,i)\mapsto \sigma .i
 =\sigma(i)$. Then our objective in below is to establish the
 spherical transform of the pair ($\Snp,\Sn$). Indeed, we will
 characterize  the right cosets $\Snp \slash \Sn$ and the double
 cosets $\Sn\backslash\Snp\slash\Sn$.
\subsection{Realization of $\Snp \slash \Sn$ of right cosets}

We now study the set of
 right cosets $\Snp \slash \Sn=\{\sigma.\Sn ;\sigma \in \Snp \}$.
Consider the function $f: \Snp \to Z_{n+1}$
 defined by $f(\sigma)=\sigma(n+1)$.
As $\Snp$ acts transitively on $Z_{n+1}$
 it follows that $f$ is surjective.\\
 For $\sigma,\sigma' \in \Snp$ we have
 $f(\sigma)=f(\sigma')$ if and only if
 $\sigma^{-1}\circ\sigma'(n+1)=n+1$,
 or equivalently if $\sigma\Sn= \sigma'\Sn$.
Thus $f$ induces a bijection
 $\overline{f}:\Snp \slash \Sn\to Z_{n+1}$ given by
 $\overline{f}(\sigma\Sn)=\sigma(n+1)$. Thus $\Snp \slash \Sn = Z_{n+1}.$
\subsection{Realization of $\Sn\backslash\Snp\slash\Sn$ of double cosets}
We now consider the set of double cosets: $$
 \Sn\backslash\Snp \slash \Sn
 =
 \{ \Sn\sigma\Sn:
 \sigma \in \Snp \}.
$$
 We shall calculate for $\sigma\in\Snp$
 the double coset $\Sn\sigma\Sn$.\\
 \\
If $\sigma\in\Sn$ then $\Sn\sigma\Sn=\Sn$
 so we assume $\sigma\notin\Sn$.
Thus $\sigma(n+1)\ne n+1$ so we must have
 $\sigma(n+1)\in Z_n$.
As $\Sn$ acts transitively on $Z_n$ there is a $\sigma'\in\Sn$
 such that $\sigma'\sigma(n+1)=1$.
This means $\sigma'\sigma(n+1)=\tau_{1,n+1}(n+1)$,
 so by the discussion above we have
$$
 \sigma'\sigma \Sn = \tau_{1,n+1} \Sn.
$$ This implies $$
 \Sn\sigma'\sigma \Sn = \Sn\tau_{1,n+1} \Sn.
$$ However since $\sigma'\in \Sn$,
 we have
$$
 \Sn\sigma \Sn = \Sn\tau_{1,n+1} \Sn.
$$ Therefore there are only two double cosets. The result follows
since a group may always be
 expressed as the disjoint union of its
 double cosets with respect to any subgroup. Then
$$
 \Snp =\Sn \bigcup \Sn.\tau_{1,n+1}.\Sn,
$$
 and the union is disjoint.\\
 \\
We therefore have
 $\Sn\backslash \Snp \slash\Sn=\{\Sn ,\Sn.\tau_{1,n+1}.\Sn\}$
 so the space of radial and integrable function
 $\bL^{1,\#}(\Snp )$
 is equal to $\bL^{1}(\Sn\backslash\Snp \slash \Sn)=\bL^{1}(\{\Sn ,\Sn.\tau_{1,n+1}.\Sn\})$.\\
\\
We say that ($\Snp ,\Sn)$ is a Gelfand pair
 when the convolution algebra
 $\bL^{1,\#}(\Snp )=\bL^{1}(\Sn\backslash \Snp \slash \Sn)$
 of integrable and $\Sn$-biinvariante functions on
 $\Snp$ is abelian.\\
 \begin{Rem}
{\bf 1)} As $\Sn \backslash \Snp \slash \Sn $ is finite, it is
easy to see
 that the convolution algebra
 $\bL^{1,\#}(\Snp )=\bL^{1}(\Sn\backslash \Snp \slash \Sn)$ is
 Abelian, thus ($\Snp ,\Sn)$ is a Gelfand pair.\\
{\bf 2)} As the cardinal of $\Sn \backslash \Snp \slash \Sn$ is
equal to number of $\Snp$-orbits in $\Sn\backslash\Snp \times \Snp
\slash \Sn $ ( $\Snp$ acts via $(\Sn \sigma , \tau \Sn)
\longrightarrow (\Sn \sigma \zeta ,\zeta^{-1} \tau \Sn)$), then we
have two $\Snp$-orbits in $\Sn\backslash\Snp \times \Snp \slash
\Sn $
 \end{Rem}
 For any subset $A\subset \Snp$, we define $\chi_{A}$ be the
 characteristic function of $A$
 \begin{Cor} The element
$(n+1)\chi_{\Sn}$ is an identity element of $\bL^{1,\#}(\Snp)$.
 \end{Cor}
\begin{proof} Choose any $f\in \bL^{1,\#}(\Snp)$, we must show that
 $\chi_{\Sn}*f=f*\chi_{\Sn}=\frac{1}{n+1}f$.
By definition of convolutions we have $$
 (\chi_{\Sn}*f)(x)
 =
 \int_{\Snp}
 \chi_{\Sn}(y)f(xy^{-1})d\mu(y)
 =
 \int_{\Sn} f(xy^{-1})d\mu(y).
$$ As $f$ is right-invariant we have $$
 (\chi_{\Sn}*f)(x)
 =
 \int_{\Sn} f(x)d\mu(y)
 =
 f(x)\mu(\Sn).
$$ On the other hand $\mu(\Sn)=[\Snp:\Sn]^{-1}=\frac{1}{n+1}$,
 so we have $\chi_{\Sn}*f=\frac{1}{n+1}f$.
The formula $\chi_{\Sn}*f=\frac{1}{n+1}f$ follows in the
 same way but using left-invariance rather than right-invariance of
 $f$.
\end{proof}\\
For any subset $A\subseteq \Snp$
 we define $\chi_A$ to be the characteristic function of $A$.
We denote by $\chi_{\sigma}$ the characteristic function of
 the set $\{\sigma\}$ and by $\delta_{\sigma}$ the Dirac measure at the point
 $\sigma$.
Note that for $\sigma ,\tau \in\Snp$ we have
 $\delta_{\sigma} * \delta_{\tau}= \delta_{\sigma \tau}$.
From our normalization of the Haar measure $dk$ on $\Snp$
 it follows that $\chi_{\sigma}(dk)=\frac{1}{(n+1)!}\delta_{\sigma}$.
We therefore have $\chi_{\sigma} * \chi_{\tau}= \frac{1}{(n+1)!}
\chi_{\sigma \tau}$.\\
\\
Let $\sigma \in \Snp$, we have
\begin{eqnarray*}
 \chi_{\sigma}^{\#}(x)
 &=&
 \int_{\Sn}\int_{\Sn}
 \chi_{\sigma}(\tau.x.h)\,d\tau dh\\
 &=&
 \chi_{\Sn.\sigma.\Sn}(x).
\end{eqnarray*}
So, for $\sigma = id_{\Sn}$ $$
 \chi^{\#}_{id(\Sn)}=\chi_{\Sn}
$$ By virtue of
 $\tau_{1,n+1}=\tau_{1,i}\circ \tau_{i,n+1}\circ \tau_{i,1}$ for all $i\leq n$,
 we will have
$$
 \chi^{\#}_{1,n+1}
 =
 \chi^{\#}_{i,n+1}=\chi_{\Sn.\tau_{1,n+1}.\Sn}
$$
 and
\begin{eqnarray*}
 \chi^{\#}_{n+1,n+1}
 &=&
 \chi^{\#}_{id(\Sn)}\\
 &=&
 \chi_{\Sn}.
\end{eqnarray*}
\begin{Rem}
From Corollary 3.1, we have $$
 \chi^{\#}_{1,n+1}*\chi^{\#}_{id(\Sn)}(x)=
 \frac{1}{n+1} \chi^{\#}_{1,n+1}(x)=
 \chi^{\#}_{id(\Sn)}*\chi^{\#}_{1,n+1}(x),
$$ so any two basis elements of the convolutions algebra
$\bL^{1,\#}(\Snp)$ commute, which is another way to see that
$(\Sn,\Snp)$ is a Gelfand Pair. Also we have $$
 \chi^{\#}_{id(\Sn)}*\chi^{\#}_{id(\Sn)}= \frac{1}{n+1}\chi_{id(\Sn)}. $$
\end{Rem}
\begin{Lem}
We have $$
 \chi^{\#}_{1,n+1}*\chi^{\#}_{1,n+1}
 =
 \frac{n-1}{n+1}\chi^{\#}_{1,n+1}
 +
 \frac{n}{n+1}\chi^{\#}_{id(\Sn)}
$$
\end{Lem}
\begin{proof}
To prove this lemma, we use the fact that
 $\Snp =\Sn \bigcup \Sn.\tau_{1,n+1}.\Sn$ is
 disjoint union.
Then the characteristic function on $\Snp $
 noted by $\I_{\Snp }$ may be expressed as
$$
 \I_{\Snp }
 =
 \chi^{\#}_{1,n+1}
 +
 \chi^{\#}_{id(\Sn)}.
$$ We use the fact that $$
 \I_{\Snp } * \I_{\Snp }
 =
 \I_{\Snp }.
$$ Therefore $$
 (\chi^{\#}_{1,n+1}+ \chi^{\#}_{id(\Sn)}) *
 (\chi^{\#}_{1,n+1}+ \chi^{\#}_{id(\Sn)})
 =
 ( \chi^{\#}_{1,n+1}+ \chi^{\#}_{id(\Sn)}).
$$ Expanding the left hand side of the above equality we obtain:
$$
 \chi^{\#}_{1,n+1}*\chi^{\#}_{1,n+1}
 +
 2(\chi^{\#}_{1,n+1}*\chi^{\#}_{id(\Sn)})
 +
 \chi^{\#}_{id(\Sn)}*\chi^{\#}_{id(\Sn)}.
$$ This is equal to: $$
 \chi^{\#}_{1,n+1}*\chi^{\#}_{1,n+1}
 +
 \frac{2}{n+1}\chi^{\#}_{1,n+1}
 +
 \frac{1}{n+1}\chi^{\#}_{id(\Sn)}.
$$ We therefore have: $$
 \chi^{\#}_{1,n+1}
 +
 \chi^{\#}_{id(\Sn)}
 =
 \chi^{\#}_{1,n+1}*\chi^{\#}_{1,n+1}
 +
 \frac{2}{n+1}\chi^{\#}_{1,n+1}
 +
 \frac{1}{n+1}\chi^{\#}_{id(\Sn)}.
$$ Consequently $$
 \chi^{\#}_{1,n+1}*\chi^{\#}_{1,n+1}
 =
 \frac{n-1}{n+1}\chi^{\#}_{1,n+1}
 +
 \frac{n}{n+1}\chi^{\#}_{id(\Sn)}
$$
\end{proof}

\subsection{Spherical function of the Gelfand pair $(\Snp ,\Sn)$}

A function $\phi$ is said to be a spherical function (see [7])
 if and only if $\phi$ is $\Sn$-biinvariante
 and $\phi$ is a character of $\bL^{1,\#}(\Snp )$.
Then a function $\phi$ of $\Snp $ which is
 $\Sn$-biinvariante may be considered as a function of\\
 $\Sn\backslash \Snp \slash \Sn=  \{\Sn , \Sn\tau_{1,n+1}\Sn\}$.
We shall use the following notation
 $\check{f} (x) = f(x^{-1})$ and
 $\widetilde{(f)}(x)=\overline{f(x^{-1})}.$
\begin{Teo}
The spherical functions of the Gelfand pair
 $(\Snp ,\Sn)$ are of the form
\begin{itemize}
 \item[1)]
 The characteristic function $\I=\chi_{1,n+1}^{\#}+\chi^{\#}_{id(\Sn)}$
  on $\Snp$ whose restriction
  to $\Sn$ is equal to $\chi_{id(\Sn)}$;
 \item[2)]
 The function
  $\phi_{n}=\frac{-1}{n}\chi_{1,n+1}^{\#}+\chi^{\#}_{id(\Sn)}$
\end{itemize}
\end{Teo}
\begin{proof}
Let $\phi$ be a spherical function, then $\phi$ is an
 $\Sn$-biinvariante function such that $\phi(Id_{\Snp})=1$ and
 satisfying the following integral equation (see [7])
$$
 \phi(\sigma)\phi(\zeta)
 =
 \int_{\Sn}
 \phi(\sigma.\tau.\zeta)d\nu(\tau),
$$ where $\sigma$ and $\zeta \in \Snp$. Also we have,
 $\widetilde{\phi}=\phi$ and
 $\check{f}*\phi=<f,\phi>\phi$,\,\,\,
 $\forall f \in \bL^{1,\#}(\Snp )$, with $<f,g>=\frac{1}{(n+1!)}\Sigma_{\sigma \in \Snp}f(\sigma)g(\sigma)$.\\
As $\phi$ is biinvariant we may express it in terms of our
 basis:
$$
 \phi
 =
 \alpha \chi_{1,n+1}^{\#}
 +
 \beta \chi^{\#}_{id(\Sn)}.
$$ As $\phi(e)=1$ we must have $\beta=1$. Therefore $$
 \phi
 =
 \alpha \chi_{1,n+1}^{\#}
 +
 \chi^{\#}_{id(\Sn)}.
$$ As $\check{\phi} * \phi = \overline{\phi}* \phi = <\phi , \phi>
\phi$,
 we have, in the first hand
\begin{eqnarray*}
 \check{\phi} * \phi
 &=&
 (\overline{\alpha} \chi_{1,n+1}^{\#}
 +
 \chi^{\#}_{id(\Sn)})* (\alpha \chi_{1,n+1}^{\#}
 +
 \chi^{\#}_{id(\Sn)})\\
 &=&
 |\alpha |^{2}(\chi_{1,n+1}^{\#} + \chi^{\#}_{id(\Sn)})
 +
 2 Re(\alpha)(\chi_{1,n+1}^{\#} + \chi^{\#}_{id(\Sn)})
 +
 (\chi_{1,n+1}^{\#} + \chi^{\#}_{id(\Sn)})\\
 &=&
 |\alpha |^{2}(\frac{n-1}{1+n}\chi_{1,n+1}^{\#}
 +
 \frac{n}{1+n}\chi^{\#}_{id(\Sn)})
 +
 2Re(\alpha)\frac{1}{1+n}\chi_{1,n+1}^{\#}
 +
 \frac{1}{1+n}\chi^{\#}_{id(\Sn)}\\
 &=&
 \frac{(n-1)|\alpha |^{2} +2 Re(\alpha)}{n+1}\chi_{1,n+1}^{\#}
 +
 (\frac{n|\alpha |^2 +1}{n+1})\chi^{\#}_{id(\Sn)}.
\end{eqnarray*}
In the second hand $$
 <\phi, \phi>
 =
 <\alpha \chi_{1,n+1}^{\#} + \chi^{\#}_{id(\Sn)}, \alpha \chi_{1,n+1}^{\#} +
\chi^{\#}_{id(\Sn)}>. $$ So
\begin{eqnarray*}
 <\phi, \phi>
 &=&
 |\alpha |^2<\chi_{1,n+1}^{\#},\chi_{1,n+1}^{\#}>
 +
 <\chi^{\#}_{id(\Sn)},\chi^{\#}_{id(\Sn)}>\\
 &=&
 \frac{|\alpha |^2}{(n+1)!}
 \sum_{\sigma \in \Snp }
 \chi_{1,n+1}^{\#}(\sigma )\chi_{1,n+1}^{\#}(\sigma )\\
 &+&
 \frac{1}{(n+1)!}\sum_{\sigma \in \Snp}
 \chi^{\#}_{id(\Sn)}(\sigma)\chi^{\#}_{id(\Sn)}(\sigma)\\
 &=&
 \frac{|\alpha |^2 n n!}{(n+1)!} +  \frac{n!}{(n+1)!}\\
 &=&
 \frac{n |\alpha |^2 +1}{(n+1)}.
\end{eqnarray*}
Then
\begin{eqnarray*}
 <\phi, \phi> \phi
 =
 (\frac{n |\alpha |^2 +1}{(n+1)}) \alpha \chi_{1,n+1}^{\#}
 +
 (\frac{n |\alpha |^2 +1}{(n+1)})\chi^{\#}_{id(\Sn)}.
\end{eqnarray*}
By vertue of $\check{\phi}* \phi = <\phi,\phi>\phi $, we will have
\begin{eqnarray*}
 (\frac{(n-1)|\alpha |^2+ 2n Re(\alpha)}{n+1})
 &=&
 (\frac{n |\alpha |^2+1}{n+)})\alpha \,\,\,\,and\\
 (\frac{n |\alpha |^2 +1}{n+1})
 &=&
 (\frac{n |\alpha |^2 +1}{n+1})
\end{eqnarray*}
Then $(n-1)|\alpha |^{2} +2Re(\alpha)=(n|\alpha|^2+ 1)\alpha$.
From this equality, we deduce that $\alpha$ is real
 and $(n-1)\alpha^{2} +2\alpha=(n\alpha ^2+ 1)\alpha$.
Thus $\alpha (n\alpha^2 + (n-1)\alpha +1-2n)=0$
 and this equality becomes
 $\alpha(\alpha -1)(\alpha +\frac{1}{n})=0$.
The solutions of this are $\alpha =0$, $\alpha =1$
 and $\alpha = -\frac{1}{n}$.
This proves the theorem
\end{proof}
\begin{Rem}
The spherical functions on finite group are all of positive type
cf.[5, p.66]
\end{Rem}
\section{The spherical Fourier transform on permutation group}

In this section, we consider the horocyclic Radon transform which
turns functions defined on  $\Snp$ into functions defined on the
set of the horocycles.\\ We note that the limit inductive $\ds
lim_{\longrightarrow}\Sn =lim_{n\geq 1}\Sn$ =$\{\sigma$ bijection
from $I\!\!N^* \to I\!\!N^*$ such that $supp(\sigma)$ is finite
$\}$ ( where $supp(\sigma)$=$\{$k such that$ \sigma(k)\neq k\}$),
so $lim_{n\geq 1}\Sn$ =$\{\sigma$ bijection from $I\!\!N^* \to
I\!\!N^*; \exists n_{\sigma}$ such that $\sigma(k)=k, \forall
k\geq n_{\sigma}+1\}$.\\ The construction of a function from its
Radon transform is a central point of study of this section and
the short way to invert the Radon transform is to note the
connection with the Fourier transform and the horocyclic Radon
transform is defind by the following formula (see [12])
$$Rf={\mathcal{F}}_{1}^{-1}\circ \widetilde{f},$$ where
${\mathcal{F}}_{1}^{-1}$ is the inverse of the finite Fourier
transform  and $\widetilde{f}$ is the spherical Fourier transform
of the Gelfand pair $(\Snp, \Sn)$. Of course we must place
hypotheses on $f$ so that the Fourier inversion formula is valid
in order to obtain the following result $$f=R^{-1}\circ
{\mathcal{F}}_{1}^{-1}\circ \widetilde{f}.$$ More explicitly:
 The spherical Fourier transform is
\begin{eqnarray*}
\wedge{f}(n)=\widetilde{f}(n)&=&\int_{\Snp}\phi_n(\sigma)f(\sigma)d\mu(\sigma)\\
&=&\frac{1}{(n+1)!}\sum_{\sigma \in
\Snp}\phi_{n}(\sigma)f(\sigma)\\ &=& \frac{1}{(n+1)!}\sum_{\sigma
\in
\Snp}f(\sigma)(\frac{-1}{n}\chi_{1,n+1}^{\#}+\chi^{\#}_{id(\Sn)})\\
&=&\frac{1}{(n+1)!}\sum_{\sigma \in \Sn}f(\sigma) -
\frac{1}{n(n+1)!}\sum_{\sigma \in \Sn
\tau_{1,{n+1}}\Sn}f(\sigma)\\
&=&\frac{1}{(n+1)}[\frac{1}{n!}\sum_{\sigma \in \Sn}f(\sigma) -
\frac{1}{n n!}\sum_{\sigma \in \Sn \tau_{1,{n+1}}\Sn}f(\sigma)].
\end{eqnarray*}
And we set
\begin{eqnarray*}
\wedge_{1}{f}(n)=\widetilde{f_1}(n)&=&\frac{1}{n!}\sum_{\sigma \in
\Sn}f(\sigma)=f*\nu(Id)\\ \wedge_{2}{f}(n) =\widetilde{f_2}(n)&=&
\frac{1}{n n!}\sum_{\sigma \in \Sn \tau_{1,{n+1}}\Sn}f(\sigma).
\end{eqnarray*}
We note that the values of two functions $\wedge_1{f}(n)$ and
$\wedge_2{f}(n)$ are thus equal to the averages of the function
$f$ over  $\Sn$ respectively over $\Sn \tau_{n+1} \Sn$: elements
of the double cosets $\Sn\backslash\Snp\slash \Sn$.
\section{Inversion formula}
The solution of reconstruction function from its spherical Fourier
transform can be written in a simple iterative form which is
computationally very tractable. Notice that the identity
$Id=Id_{S_{i}}$ (for all $1\leq i \leq n$) . Consider first how to
find $f(Id)$ from $\widetilde{f}$.
\begin{eqnarray*}
\widetilde{f}(1)&=&\frac{1}{2}f(Id)- \frac{1}{2}f(\tau_{1,2})\\
\widetilde{f}(2)&=&\frac{1}{3!}(f(Id)+f(\tau_{1,2}))-
\frac{1}{2.3!}[f(\tau_{1,3})+f(\tau_{1,3}\circ \tau_{1,2})+
f(\tau_{1,2}\circ \tau_{1,3})\\ & &+f(\tau_{1,2}\circ \tau_{1,3}
\circ \tau_{1,2})]\\ &=& \frac{1}{3!}(f(Id)+f(\tau_{1,2}))-
\frac{1}{2.3!}[f(\tau_{1,3})+f(\tau_{1,3}\circ \tau_{1,2})+
f(\tau_{1,2}\circ \tau_{1,3})\\ & &+f(\tau_{2,3})]\\
\widetilde{f}(3)&=&\frac{1}{4!}[f(Id)+ f(\tau_{1,2})+
f(\tau_{1,3})+ f(\tau_{2,3})+f(\tau_{1,3}\circ
\tau_{1,2})+f(\tau_{1,2}\circ \tau_{1,3})]\\ & &-
\frac{1}{3.4!}[f(\tau_{1,4})+f(\tau_{1,4}\circ
\tau_{1,2})+f(\tau_{1,4}\circ \tau_{2,3})+f(\tau_{1,4}\circ
\tau_{1,3})\\ & &+f(\tau_{1,4}\circ \tau_{1,3}\circ \tau
_{1,2})+f(\tau_{1,4}\circ \tau_{1,2}\circ \tau
_{1,3})+\sum_{\sigma \in S_3 \tau_{1,4} S_3}f(\sigma)]\\
......&=&.....\\ \widetilde{f}(n)&=&\frac{1}{(n+1)!}[f(Id)+
\sum_{\sigma \in \Sn \backslash {Id}}f(\sigma)]-
\frac{1}{n.(n+1)!}\sum_{\sigma \in \Sn \tau_{1,2}\Sn}f(\sigma).
\end{eqnarray*}
Since $f(Id)$ occurs in the sum for $\widetilde{f}(i)$ (for all
$1\leq i$) and $f(\tau_{1,2})$ occurs also in the sum for
$\widetilde{f}(i)$ (for all $1\leq i$) and continuing in this way
in a simple iterative form , we have an explicit inversion formula
for $f(Id)$ which can easily must be
$$f(Id)=\frac{1}{n}\sum_{k=1}^{k=n}(1+k)!\widetilde{f}(k)-\frac{1}{n}
\sum_{k=2}^{k=n}(n-k)f(\tau_{1,k})-\frac{1}{n^2}f(\tau_{1,n+1})+...$$
Applying the same reasoning for any $\tau_{1,i} \in \Snp$ in stead
of $Id$ leads to the full inversion formula
$$f(\tau_{1,i})=\frac{1}{n-i}\sum_{k=1}^{k=n}(1+k)!\widetilde{f}(k)-\frac{n}{n-i}f(Id)+...$$
\\begin{Rem}
In the case of permutations group, the difficulty arises for an
excellent account of analysis and modeling of front propagation in
permutation group and their utility in applications to the heat
equations.\section{ Horocyclic Radon Transform and inversion
formula} Given a function $f: \Snp \to \C$, the Radon transform of
$f$ is
\begin{eqnarray*}
R_{n+1}f\left(\sigma\right)&=&\sum_{\sigma' \in \sigma
\Sn}f(\sigma')\\ &=&\frac{1}{n!}\sum_{h \in \Sn}f(\sigma.h)
\end{eqnarray*}
\begin{Rem}   For any $h_0 \in \Sn$, we have
$$R_{n+1}f(\sigma)=R_{n+1}nf(\sigma.h_0),\,\,\,\,\,\,
R_{n+1}f(Id_{\Sn})=\frac{1}{n!}\sum_{h\in \Sn}f(\sigma.h)
=\wedge_1f(n)$$
\end{Rem}
 This Radon transform is also defined for every $B\in \Snp\slash \Sn=Z_{n+1}$ and
every $f\in L( Z_{n+1})$ by $$Rf(B)=\sum_{A\in Z_{n+1}: A\subset
B} f(A).$$  So for an absolutely summable function on  $\ds
lim_{\rightarrow}Z_{n+1}=\N^*$, the Radon transform (see
[16])$$Rf(m)=\sum_{k=1}^{k=\infty}f(k m),$$ for all $m\in \N^*$,
assuming some more rapid decay (say $|f(n)|<cn^{-2-\epsilon}$).\\
This transform can be easily inverted (see [16]), First we find
$f(1)$ from $Rf$ of the form $$f(1)= c_1Rf(1)+
c_2Rf(2)+...+c_nRf(n)+...$$ for certain coefficients $c_n$ that
are uniquely determined. The coefficients must be equal to the
Mobius function $\mu(n)$, which is defined to be $(-1)^k$ if $n$
has $k$ distinct prime factors, and 0 if $n$ is divisible by a
square of a prime. The full inversion formula (see [16])is
$$f(n)=\sum_{k=1}^{\infty}\mu(k)Rf(nk)$$
\section{Applications}
In order to state the applications of the Radon transform, we
consider functions $f: \Snp \times I\!\!N \longrightarrow
1\!\!\!C$,
 whose values are denoted by $f(x, k)=f_{k}(x)$. In our setting,
 $x\in \Snp$ represents the space variable and $k\in I\!\!N$ the
 times variables. Let the sub-Laplacian of $f$ with respect to the space
 variable defined in the following way
 $$\bigtriangleup f(x)=\frac{1}{n!}\sum_{h\in
 \Sn}f(x.h)-f(x)=f*\nu(x)-f(x).$$Now we consider the following boundary value
 problems, denoted ${\cal{H}}$, in
 analogy with the corresponding problems for the classical heat
equation.\\
 \begin{eqnarray*}
 ({\cal{H}}) \,\,\,\,\,\,\,\,\,\,\,\,\,\,\, \bigtriangleup
 f_{k}(x)&=&f_{k+1}(x)- f_{k}(x)
 \,\,\,\,\,\,\,\,\,\,\,\,\,\,\,given\,\,\,\,
 f_0(x).
\end{eqnarray*}
The heat equation associated to sub-Lapcian operator can be
written in the following way
$$({\cal{H}})\,\,\,\,f_{k}(x)=\frac{1}{n!}\sum_{h\in
 \Sn}f_{k-1}(x.h)\,\,\,\, given\,\,\,\, f_0=\delta_{e},$$
\begin{Teo}
The solution of the heat equation is given by
$$f_{k}(x)=f_0*\nu(x),$$
 with $\nu=\frac{1}{n!}\sum_{\sigma \in \Sn} \delta_{\sigma}$ the
 Haar measure in $\Sn$.\\If $f_0=\delta_e$, we have $f_k=\nu$
 \end{Teo}
 \begin{proof}
The equation $({\cal{H}})$ is therefore
 $$f_{k}(x)=f_{k-1}*\nu(x)=f_0*\nu^{k}(x),$$
 with $\nu=\frac{1}{n!}\sum_{\sigma \in \Sn} \delta_{\sigma}$ the
 Haar measure in $\Sn$ and $\nu^{k}=\nu*\nu*\nu*...*\nu$
 $k$-times. Using the fact that $$(\sum_{\sigma \in \Sn} \delta_{\sigma})*(\sum_{\sigma \in \Sn}
 \delta_{\sigma})=n!\sum_{\sigma \in
 \Sn} \delta_{\sigma}.$$Because $(\sum_{\sigma \in \Sn} a_\sigma \delta_{\sigma})*(\sum_{\sigma \in \Sn} b_\sigma \delta_{\sigma})=\sum_{\sigma \in
 \Sn}(\sum_{tu=\sigma}a_t b_u)
 \delta_{\sigma}.$\\
 So $(\sum_{\sigma \in \Sn} \delta_{\sigma})^k=(n!)^{k-1}\sum_{\sigma \in
 \Sn}
 \delta_{\sigma}.$\\ Then $\nu^k=\nu$ and we will have
 $f_k=f_0*\nu=\nu$, with $f_0=\delta_e$
 \end{proof}

\end{document}